\documentclass{birkjour}
\usepackage{t1enc}
\usepackage{epsf}
\usepackage{graphicx}
\usepackage[pdftex]{color}
\usepackage{amsthm}
\usepackage{amsmath,amssymb}

\renewcommand{\figurename}{{\bf fig.}}
\newcommand{\N}{\mbox{I\kern-2ptN}}
\newcommand{\R}{\mbox{I\kern-2ptR}}
\newcommand{\C}{\mbox{I\kern-6ptC}}

\newtheorem{teo}{Theorem}

\theoremstyle{definition}

\definecolor{morado}{rgb}{0.064,0,0.153}
\definecolor{azulclaro}{rgb}{0.8,0.85,1}

\begin{document}

\title{The NIEP}  

% CRJ
\author[Johnson]{Charles R.~Johnson} 
\address{
Department of Mathematics \\ 
College of William and Mary \\ 
Williamsburg, VA 23187-8795 \\ 
USA
}
\email{crjohn@wm.edu}

% CM
\author[Mariju\'an]{Carlos Mariju\'an}
\address{
Dpto. Matem\'atica Aplicada \\ 
E.I. Inform\'atica\\ 
Paseo de Bel\'en 15	\\ 
47011-Valladolid \\
Spain
}
\email{marijuan@mat.uva.es}

% PP
\author[Paparella]{Pietro Paparella}
\address{
Div. of Engineering \& Mathematics \\ 
University of Washington Bothell \\ 
Bothell, WA 98011-8246 	\\
USA}
\email{pietrop@uw.edu}

% MP
\author[Pisonero]{Miriam Pisonero}
\address{
Dpto. Matem\'atica Aplicada \\ 
E.T.S. de Arquitectura  \\
Avenida de Salamanca 18 \\ 
47014-Valladolid \\
Spain
}
\email{mpisoner@maf.uva.es}

\begin{abstract}
The nonnegative inverse eigenvalue problem (NIEP) asks which lists of $n$ complex numbers (counting multiplicity) occur as the eigenvalues of some $n$-by-$n$ entry-wise nonnegative matrix. The NIEP has a long history and is a known hard (perhaps the hardest in matrix analysis?) and sought after problem. Thus, there are many subproblems and relevant results in a variety of directions. We survey most work on the problem and its several variants, with an emphasis on recent results, and include 130 references. The survey is divided into: a) the single eigenvalue problems; b) necessary conditions; c) low dimensional results; d) sufficient conditions; e) appending 0's to achieve realizability; f) the graph NIEP's; g) Perron similarities; and h) the relevance of Jordan structure.
\end{abstract}

\thanks{Partially supported by  MTM2015-365764-C-1-P(MINECO/FEDER)}

\maketitle

%----------------------------------------------------------------------------------------------------------------------------------------------------------------------------------------------------------------------------
\section{Introduction} \label{Se1}
The { \bf N}onnegative { \bf I}nverse { \bf E}igenvalue {\bf P}roblem (NIEP) asks which collections of $n$ complex numbers (counting multiplicities) occur as the eigenvalues of an $n$-by-$n$ matrix, all of whose entries are nonnegative real numbers. This is a long standing problem that is very difficult and, perhaps, the most prominent problem in matrix analysis. Unlike many other inverse eigenvalue problems, tools are limited, and seemingly small results are very welcome. Nonetheless, the problem is very attractive and has been attacked by many excellent researchers (see references). There have been prior surveys of work on the NIEP (\cite{Sw, Mi3, BeP, ELN, Bo2, Cr, So4}), but none recently. With a great deal of important recent activity, a new survey, emphasizing this activity, will be welcome to anyone considering inquiry in this area. Our intent here is to be broad and informal.

Another intriguing aspect of the NIEP is how many avenues can be taken to gain insight into it. Thought about it spawns endless challenging and worthy specific questions. An example is how many interesting variations there are. The real NIEP (R-NIEP) restricts the question to real spectra, and the symmetric NIEP (S-NIEP) further restricts to real spectra that are realizable by symmetric nonnegative matrices. It is known \cite{J2} that the row stochastic NIEP (restriction to nonnegative matrices with row sums 1) is equivalent to the general NIEP, but the doubly stochastic NIEP (DS-NIEP), restriction to matrices with both row and column sums 1, is properly more restrictive, but no less difficult. Of course, there are variations on the DS-NIEP, as well: the real DS-NIEP and the symmetric DS-NIEP. The Jordan NIEP (J-NIEP) asks about possible Jordan canonical forms, when there are repeated eigenvalues, and the diagonalizable NIEP (D-NIEP) is the special case in which the realizing matrix is diagonalizable. Finally, there are graph NIEP's when we consider only nonnegative matrices with a given (directed or undirected) graph $G$ (G-NIEP), or only matrices subordinate to same 0-pattern. Further particularizations may also be imagined, but these are the primary ones, so far. We should also mention that the NIEP and each variation also has a trace 0 version. This considers only spectra, the sum of whose elements is 0, or, equivalently, nonnegative matrices, each of whose diagonal entries is 0. In some cases, this is sufficiently restrictive, so as to make the problem easier.

We denote by NIEP$_n$ the set of all spectra enjoyed by some $n$-by-$n$ nonnegative matrix. Such a spectrum is called {\bf realizable}  and a nonnegative matrix with the spectrum is called a {\bf realization} or {\bf realizing matrix}. The notation and terminology are similar for NIEP variants. A solution to the NIEP (or variants), for a particular $n$, is an "explicit" description of NIEP$_n$, viewed as a subset of vectors in $\C^n$. NIEP$_n$ is a closed set that is connected, and even star-shaped from the origin (or from the all 1's vector $e$). The set is also semi-algebraic, but is not generally convex. If $\Lambda=\{\lambda_1,\dots,\lambda_n\}$ is a proposed spectrum, we denote by $\mbox{Tr}(\Lambda)$ the sum of its components, counting multiplies. Other power sums $s_k(\Lambda)=\sum_{i=1}^n\lambda_i^k$, called {\bf $k$-th moments}, are also important, and we usually reserve $\lambda_1$ as the spectral radius of $\Lambda$.

The Perron-Frobenius theory of nonnegative matrices ({\it e.g.} \cite{HoJo}) provides several simple, but important, necessary conditions for the NIEP (and variations). These basic necessary conditions are presented in Section \ref{Se3}, followed by several more subtle necessary conditions. The ``JLL'' conditions are now essential to work on the problem.

The NIEP has a long history since its proposal by Kolmogorov \cite{Ko}, and there are many substantial results. We begin at the natural starting point of the single eigenvalue problem (Section \ref{Se2}) for both the row stochastic NIEP and the DS-NIEP: which individual complex numbers occur in the spectra of row stochastic and of doubly stochastic matrices. The former is solved, but there have been important advances in the description of the solution and in realizations. The latter is unsolved, beyond $n=4$, though there have been important recent developments.

Section \ref{Se4} summarizes low-dimensional complete results for certain NIEP variants. General sufficient conditions for the realizability of spectra are given in section \ref{Se5}. It has long been known that $n$-fold spectra that meet simple necessary conditions, but are not realizable, may be made realizable by appending of 0 eigenvalues. Information about this phenomenon is given in section \ref{Se6}. We turn to what is known about the graph-NIEP's in section \ref{Se7}. The new idea of Perron similarities - studying the diagonalizable NIEP's via the diagonalizing similarities - is discussed in section \ref{Se8}, and the role of Jordan structure in the NIEP and R-NIEP in section \ref{Se9}. 
%In section \ref{Se10}, we close by raising some new and recounting some old questions related to several NIEP's.

%----------------------------------------------------------------------------------------------------------------------------------------------------------------------------------------------------------------------------
\section{The Single Eigenvalue Problems} \label{Se2}
%Row stochastic (D-D-K, history, new approaches, J-K-S variants, etc), doubly stochastic, Perfect-Mirsky, recent results.

In \cite{Ko}, Kolmogorov posed the problem, denoted by SISEP ({\bf S}tochastic {\bf I}nverse {\bf S}ingle {\bf E}igenvalue {\bf P}roblem), of characterizing the subset of the complex plane, denoted by $\Theta_n$, consisting of the individual eigenvalues of all $n$-by-$n$ stochastic matrices. 

It can easily be verified that for each $n \geq 2$, the region $\Theta_n$ is closed, inscribed in the unit-disc, star-convex (with star-centers at zero and one), and symmetric with respect to the real-axis. Furthermore, it is clear that $\Theta_n \subseteq \Theta_{n+1}$, $\forall n \in \N$. In view of these properties, each region $\Theta_n$ is determined by its boundary,  which consists of so-called {\bf extremal} numbers, {\it i.e.}, $\partial \Theta_n = \{ \lambda \in \Theta_n : \alpha \lambda \not \in \Theta_n,\forall \alpha > 1\}$.  

%Dear All, It appears that we have omitted a reference that I think is important. In fact this reference should be mentioned in section 2, where [58] is referenced (more later): C. R. Johnson, R. B. Kellogg, and A. B. Stephens, Complex Eigenvalues of a Nonnegative matrix, with a Specified Graph II, Linear and Multilinear Algebra 7 (1979), 129-143. This is a follow-on to [58] and also related to [47]. Where [58] is referenced in Sec 2, this paper should also be referenced. It has more extensive information about such things. In the sentence in section 2 where [58] is referenced and there is mention of the longest circuit of G, we should probably say "simple" circuit to avoid any confusion.

In pursuit of characterizing $\Theta_n$, Dmitriev and Dynkin \cite{DD1} (and independently Karpelevich \cite{Ka}) showed that if $\lambda = a+ bi \in \Theta_n$ and $b \neq 0$, then 
\[ a + |b| \tan\left(\frac{\pi}{n}\right) \leq 1. \] 
Let $A$ be a nonnegative matrix of order $n$ with spectral radius $\rho$ and associated directed graph $G$. Let $m$ be the length of the longest simple circuit of $G$. Kellogg and Stephens \cite[Theorem 1]{KeSt} showed that if $m = 2$, all eigenvalues of $A$ are real. If $2 < m \leq n$, and if $\lambda = a + bi$ is an
eigenvalue of $A$, then 
\[ a + |b|\tan\left(\frac{\pi}{m}\right) \leq \rho. \] 
For further results on this topic, see \cite{JKS}.

Karpelevich \cite[Theorem B]{Ka}, expanding on the work of Dmitriev and Dynkin in \cite{DD2}, resolved SISEP by showing that the boundary of $\Theta_n$ consists of curvilinear arcs (called K-arcs), whose points satisfy a polynomial equation that is determined by the endpoints of the arc (which are consecutive roots of unity). The statement of this theorem is unwieldy, but Ito \cite[Theorem 2]{Ito} provides a useful simplification of the result. 

Noticeably absent in the Karpelevich Theorem (and other works) are \emph{realizing-matrices} ({\it i.e.}, a matrix whose spectrum contains a given eigenvalue) for points on the K-arcs. Dmitriev and Dynkin \cite[Basic Theorem]{DD2} give a schematic description of such matrices for points on the boundary of $\Theta_n \backslash \Theta_{n-1}$ and Swift \cite[\S 2.2.2]{Sw} provides such matrices for $3\leq n \leq 5$. Recently, Johnson and Paparella \cite{JoPa2} provide, for every $n$ and for each K-arc, a single parametric stochastic matrix that realizes the entire K-arc as the parameter runs from 0 to 1.

For $n \in \N$, denote by $\Pi_n$ the convex-hull of the $n$\textsuperscript{th} roots-of-unity, {\it i.e.}, 
\[ \Pi_n = \left\{ \sum_{k=0}^{n-1} \alpha_k \exp{(2\pi i k/n)} \in \C : \alpha_k \geq 0,~\sum_{k=0}^{n-1} \alpha_k =1 \right\}. \]
Denote by $\Omega_n$ the subset of the complex-plane containing all single eigenvalues of all $n$-by-$n$ doubly stochastic matrices. Perfect and Mirsky \cite{PeMi} conjectured that $\Omega_n = \bigcup_{k=1}^n \Pi_k$ and proved their conjecture when $1 \leq n \leq 3$. Levick et al.~\cite{LPK} proved the Perfect-Mirsky conjecture when $n=4$ but a counterexample when $n=5$ was given by Mashreghi and Rivard \cite{MaRi}. Recently, Levick {\it et al.} conjectured that $\Omega_n = \Theta_{n-1} \cup \Pi_n$ (\cite[Conjecture 1]{LPK}), but there is computational evidence that suggests that the $n=5$ case is either a rare exception, or the only exception, to Perfect-Mirsky.

%----------------------------------------------------------------------------------------------------------------------------------------------------------------------------------------------------------------------------
\section{Necessary Conditions for the NIEP's} \label{Se3}
For a proposed spectrum $\Lambda=\{\lambda_1,\lambda_2,\dots,\lambda_n\}$, repeats allowed, to be (NIEP-)realizable, a number of necessary conditions are known. The most basic of these follow from the fact that a nonnegative matrix has real entries and nonnegative trace, and from the Perron-Frobenius theory of nonnegative matrices. (We assume the reader is familiar with the Perron-Frobenius theory, which is recounted, for example in \cite{HoJo}.) For simplicity, we label these according  to simple titles. Since a nonnegative matrix is real, it must have the eigenvalues of a real matrix, {\it  i.e.} the characteristic polynomial must have real coefficients, or, as a list

\vspace*{.3cm}
\noindent $\mbox{(Reality) } \hspace{4cm} \bar{\Lambda}=\Lambda.$
\vspace*{.3cm}

\noindent Since the trace of a matrix is the sum of the eigenvalues and the trace of a nonnegative matrix is nonnegative, we also have 

\vspace*{.3cm}
\noindent $\mbox{(Trace) } \hspace{3.5cm} \mbox{Tr }\Lambda=\sum\limits_{i=1}^n\lambda_i\geq 0.$
\vspace*{.3cm}

\noindent Of course, it follows that if $\Lambda$ is realizable and $\mbox{Tr }(\Lambda)=0$, then all diagonal entries of any realizing matrix must be 0. Higher power sums are also nonnegative as they are traces of positive integral powers of nonnegative matrices: 

\vspace*{.3cm}
\noindent $\mbox{(k-th moment)} \hspace{1.5cm} s_k(\Lambda)=\lambda_1^k+\cdots+\lambda_n^k\geq 0, \quad k\geq 1.$
\vspace*{.3cm}

According to the Perron-Frobenius theory, the spectral radius of a nonnegative matrix must, itself, be an eigenvalue, so among the eigenvalues $\Lambda$ must be a nonnegative one, at least as big in absolute value as any others. Without loss of generality, this one may be taken to be $\lambda_1$. So, we have

\vspace*{.3cm}
\noindent $\mbox{(Perron) } \hspace{2.5cm} \lambda_1\geq |\lambda_i|,\quad i=2,\dots,n,$
\vspace*{.3cm}

\noindent and if any of the other eigenvalues are the same as $\lambda_1$, any realizing matrix must be reducible. Moreover, if all the inequalities are strict, any irreducible realizing matrix must be primitive.

The above conditions are necessary for all the NIEP's, but there are some necessary conditions for particular NIEP's.

Further necessary conditions are more subtle, but many have been noticed. The most ubiquitous of these was noticed independently in \cite{J2} and in \cite{LoLo} and is usually referred to as the JLL conditions. They generalize the trace condition and follow from the fact that every power of a nonnegative matrix is nonnegative and that positive diagonal entries must contribute to positive diagonal entries in powers ({\it e.g.} $s_p(\Lambda)>0$ implies $s_{pq}(\Lambda)>0$ for  integers $p\geq 1$ and $q\geq 1$). The general (quantitative) version is 

\vspace*{.3cm}
\noindent $
\mbox{(JLL) } \hspace{2cm} (s_k(\Lambda))^m \leq n^{m-1}s_{km}(\Lambda),\quad k,\, m = 1, 2, \dots .$
\vspace*{.3cm}

Let $A$ be an  $n$-by-$n$ real matrix with spectrum  $\Lambda= \{\lambda _1, \dots,\lambda _n\}$. 
Denote the principal submatrix of $A$ lying in the rows and columns given by the index set $\alpha \subseteq \{1,\dots ,n\}$ by $ A[\alpha].\,$ 
Define the $k$-th elementary symmetric function 
$$
 E_k(\Lambda )= \displaystyle\sum_{1\leq i_1<\cdots < i_k \leq n}\lambda_{i_1}\cdots \lambda_{i_k} 
$$ 
\noindent and the $k$-th Newton coefficient
$$
c_k(\Lambda )=\displaystyle \frac{E_{k}(\Lambda)}{{n \choose k}}, \quad k=1, \cdots ,n, \;\; \mbox{with} \;\; c_0\equiv 1.
$$
\noindent Since $E_k(\Lambda)=\sum_{|\alpha |=k} \det A[\alpha ]$, as well, $c_k(\Lambda)$ may be viewed as the average value of 
the $k$-by-$k$  principal minors of $A$.
The spectrum $\Lambda$ is called {\bf Newton} if
\begin{eqnarray*}\label{Ho}
 c_k(\Lambda)^2\geq c_{k-1}(\Lambda)c_{k+1}(\Lambda),\;  \quad k=1,\dots, n - 1,
\end{eqnarray*}
\noindent and these inequalities are referred to as the {\bf Newton inequalities} \cite{Ho,JMP}.

 They hold for %sequences of real numbers, for matrices with real spectrum and for 
 $\Lambda \geq 0$ \cite{Nw}, $\Lambda \subset \R $ \cite{MacL}, and they are valid for
 real diagonal matrices, diagonalizable matrices with real spectra (the $c_k$ are invariant under similarity),
 and matrices with real spectra. 
In \cite{Ho} it was proved that the Newton inequalities also hold
 for $M$-matrices and, thus, inverse $M$-matrices and it was observed that if $A$ is a nonnegative matrix with spectral radius $\rho(A)$, 
then $\rho(A)I - A$ is an M-matrix and, therefore, must satisfy the Newton inequalities.
If we denote its spectrum $\{\rho (A)-\lambda _1, \dots,\rho (A)-\lambda _n\}$ by $\rho (A)  - \Lambda$, then 
we have new necessary conditions:

\vspace*{.3cm}
\noindent $
\mbox{(H) } \hspace{.3cm}  c_k(\rho (A)-\Lambda)^2\geq c_{k-1}(\rho (A)-\Lambda)c_{k+1}(\rho (A)-\Lambda),\quad  k=1,\dots, n - 1.$
\vspace*{.3cm}

%\noindent Coefficient conditions (2007) 
 
\vskip0.3cm
 In \cite{TAAMP} the authors  did not focus the attention directly on the spectrum but on
the coefficients of the characteristic polynomial. Thus, the NIEP that they consider is:
{\it given real numbers $k_1, k_2, \dots , k_n$, find necessary and sufficient 
conditions for the existence of a nonnegative matrix of order $n$  
with characteristic polynomial $x^n+k_1x^{n-1}+k_2x^{n-2}+\dots +k_n$}.

The coefficients of the characteristic polynomial are closely related to the cyclic structure of the weighted digraph associated with the matrix $A$, 
as established by the Coefficients Theorem \cite[Theorem 1.3*]{CDS}. The authors \cite{TAAMP} introduce graphic tools to study the NIEP from the characteristic polynomial and 
use the following method:  
%the coefficient $k_1$ is minus the trace of a nonnegative matrix, so $k_1$ must be non positive. Now, for a fixed nonpositive $k_1$ we can calculate 
%the maximum possible value for $k_2$. In general, 
if  $P(x)$ is a realizable polynomial, in the sense that there exists a nonnegative matrix with 
characteristic polynomial $P(x)$, we try to maximize each coefficient $k_j$ as a function of the previous coefficients, preserving the realizability 
for a polynomial of degree $n$ with the same  previous coefficients.
Note that $k_j$ is a continuous function (sum of determinants) of the entries of a nonnegative matrix $A$ 
realizing the polynomial $P(x)$ and  that, as the previous coefficients are bounded above, then the entries of $A$ involved in the expression of  $k_j$
are also bounded above (Coefficients Theorem); therefore, this maximum is attained. In this way, new necessary conditions on the three first coefficients
are obtained \cite[Theorem 3]{TAAMP}:
\begin{eqnarray*}\label{k1}
 k_{1}\leq 0; %\addtocounter{eqnarray}{1}\hfill(\theequation).
\end{eqnarray*}

\vspace*{.3cm}
\noindent $
\mbox{(TAAMP) } \hspace{3.25cm}  k_2\leq \frac{n-1}{2n}\,k_{1}^{2};$%\addtocounter{equation}{1}\hfill(\theequation).$
\vspace*{.3cm}

 \begin{equation*}
k_3\leq \left\{\begin{array}{lll}
\frac{n-2}{n}\left(k_{1}k_{2}+ \frac{n-1}{3n}\left(\left(
k_{1}^{2}-\frac{2nk_{2}}{n-1}\right)^{\frac{3}{2}}-k_{1}^{3}\right)\right)
& {\rm if} & \frac{(n-1)(n-4)}{2(n-2)^{2}}\,k_{1}^2< k_{2},\\
\vspace*{-.3cm}& & \\
k_{1}k_{2}-\frac{(n-1)(n-3)}{3(n-2)^2}\,k_{1}^{3} & {\rm if} &
k_{2}\leq \frac{(n-1)(n-4)}{2(n-2)^{2}}\,k_{1}^2.\\
\end{array}\right.
%\addtocounter{equation}{1}\hfill(\theequation)
\end{equation*}

%\noindent Moreover, given $k_1$, $k_2$ and $k_3$ verifying the above conditions there exists a nonnegative matrix of order $n$ 
%whose characteristic polynomial is of the form $x^{n}+k_{1}x^{n-1}+k_{2}x^{n-2}+k_3x^{n-3}+Q(x)$, where $Q(x)=0$ if 
%$n=3$ and a polynomial of degree lower than or equal to $n-4$ if $n>3$.

%\noindent Cronin-Laffey condition (2012)

\vskip0.3cm
Another necessary condition in terms of the $k$-th moments $s_k=\mbox{Tr }(A^k)$, $k=1,2,3$, is obtained in \cite[Theorem 3]{CrLa}:

\vspace*{.3cm}
\noindent $
\mbox{(CL) } \hspace{1.5cm} \Phi :=n^2s_3-3ns_1s_2+2s_1^3+\frac{n-2}{\sqrt{n-1}}(ns_2-s_1^2)^{3/2}\geq 0.$
\vspace*{.3cm}
 
%\noindent {\bf Their dependency/independency}
 
The necessary conditions previous to (JLL) are not independent. In fact, ($k$-th moment) implies (Perron) \cite{Fri} and also (Reality) \cite{LoLo}.
 For $k=1,$ (JLL) is reduced to $(s_{1}(\Lambda))^{m}\leq n^{m-1}s_{m}(\Lambda)$, and so, if $s_{1}(\Lambda)\geq 0,$ then  (JLL) implies ($k$-th moment).
%Therefore, the historical necessary conditions can be reduced to (Trace) and  (JLL)  or only to (JLL) if we take $\Lambda$ with $s_1(\Lambda)\geq 0.$

On the other hand, if we denote the two bounds given for the coefficient $k_3$ in (TAAMP) by  $k_3^{max1}$ and $k_3^{max2}$, we can use the
Newton identities 
\begin{eqnarray}\label{NwI}
 s_{m}+k_{1}s_{m-1}+\dots +k_{m-1}s_{1}+k_mm=0, \quad k,m= 1,2,3,
\end{eqnarray}
to rewrite the condition (CL) in the form
\begin{eqnarray*}\label{CroLa}
\Phi =3n^2(k_3^{max1}-k_3).
\end{eqnarray*}
\noindent In the first case, {\it i.e.}, if 
${\normalsize k_2 >\frac{(n-1)(n-4)}{2(n-2)^{2}}\,k_{1}^2}$,
we have that $\Phi \geq 0$ implies $k_3^{max1}\geq k_3$. In the contrary case, $\Phi = 3n^2(k_3^{max1}-k_3)\geq  3n^2(k_3^{max2}-k_3)\geq 0 $.
So, in any case, the condition on the coefficient $k_3$ is stronger than the condition (CL).

In \cite{Ho} it was proved that ($k$-th moment), (JLL) and (H) are mutually independent.
In \cite[Theorem 11]{MP2010} that (JLL)  for  $k=1$ and $m=2$,  $s_1(\Lambda)^2\leq ns_{2}(\Lambda)$,
the first (H), $c_1(\rho-\Lambda)^2\geq c_2(\rho-\Lambda)$, and
 (TAAMP) over the second coefficient, $k_2(\Lambda)\leq \frac{n-1}{2n}k_1(\Lambda)^2$, are equivalent. But, in general, at least 
(JLL) and (TAAMP) are independent of the others, and this tandem implies ($k$-th moment). 
Some examples: the spectra $ \{20, -18, 5\sqrt{2}\pm 5\sqrt{2}i\}$ and $\{1,1,1,0,0\}$ satisfy (JLL), (H) and  (CL), but not (TAAMP);  
the spectra $ \{2, -2,-2, 1\pm i\}$ and $\{3,1,1,1,1,1,-2,-2,-2,-2\}$ satisfy (H) and  (TAAMP), but not (JLL). 
That the necessary conditions (JLL), (H) and (TAAMP) are not sufficient for the NIEP is proved by the non-realizable list
$\{3, 3, -\sqrt{3}\pm i\}$. 

\vspace*{.3cm}
\noindent {\bf Conjecture}. (JLL) and (TAAMP) imply all known necessary conditions. (It is enough to prove that (JLL) and (TAAMP) imply (H)).

%\noindent {\bf NIEP trace 0}.

\vskip0.3cm
In \cite{LM1} a necessary  condition for trace 0 and $\,n\,$ odd was obtained, given by 

\vskip0.3cm
\noindent 
$\mbox{(LM) } \hspace{3cm} (s_2(\Lambda))^2 \leq (n-1)s_{4}(\Lambda).$

\vskip0.3cm
This condition has been generalized in terms of the coefficients by the following result \cite[Lemma 37]{TAAMP}:
if $x^{n}+k_p x^{n-p}+\cdots+k_{2p} x^{n-2p}+\cdots+k_n,\; k_p\neq 0$, is the characteristic polynomial
 of a nonnegative matrix, then 
\begin{eqnarray}\label{L37}
 k_{2p}\leq \frac{1}{2}\Big(1-\frac{1}{\lfloor n/p\rfloor}\Big) k_p^2.
\end{eqnarray}
We can use the Newton  identities (\ref{NwI}) to express the inequality (\ref{L37}) in terms of the $k$-th moments: 
if $s_1=\cdots =s_{p-1}=0$, we have  the more general condition
\begin{eqnarray*}
 s_p^2 \leq p\Big\lfloor \frac{n}{p} \Big\rfloor s_{2p}, \;\; p=1,\dots, \frac{n}{2}
\end{eqnarray*}
that coincides with (LM) in the particular case $p=2$ and $n$ odd.
%   $s_2^2 \leq 2\Big\lfloor \frac{n}{2} \Big\rfloor s_{4}=(n-1)s_{4}$. 

Note also that the (JLL) condition for $m=2$ is $s_p^2 \leq ns_{2p},$ and that $p\Big\lfloor \frac{n}{p} \Big\rfloor \leq n$.
Then the expression (\ref{L37}) is a restricted refinement of the (JLL) conditions.

%----------------------------------------------------------------------------------------------------------------------------------------------------------------------------------------------------------------------------
\section{Low Dimensional Results} \label{Se4}

The NIEP for $n\leq 3$ was solved independently by Oliveira \cite[Theorem (6.2)]{O}  and Loewy-London \cite{LoLo}.
For the non-real case:
$$
\{\lambda ,z, \bar{z}\} \;\; \mbox{is NIEP-realizable}  \; \Longleftrightarrow \; z\in \lambda \Pi _3=\{\lambda z': z'\in \Pi _3\}
$$
Meehan \cite{Me} solved the NIEP for $n=4$ in terms of the $k$-th moments, 
and Torre-Mayo {\it et al.} \cite{TAAMP} in terms of the coefficients of the characteristic polynomial
using the necessary conditions (TAAMP).
For $n\geq 5$ it remains unsolved.

The R-NIEP for $n\leq 4$ was solved independently by Perfect \cite{Pe1}, Oliveira \cite[\S9] {O} and Loewy-London \cite{LoLo}; 
in these cases the necessary conditions (Trace) and (Perron) are also sufficient.

Perfect \cite[Theorem 4]{Pe3} solved the R-NIEP for $n=3$ with fixed diagonal entries:
$\{d_1,d_2,d_3\}$ is the diagonal of a 3-by-3 nonnegative matrix with spectra $\{\lambda _1, \lambda _2, \lambda _3\}$, where $\lambda _1\geq \lambda _2\geq \lambda _3$, if and only if
$$
0\leq d_i\leq \lambda _1;  \;\; \sum_{i=1}^3d_i=\sum_{i=1}^3\lambda _i; \;\; \sum_{i\neq j} d_id_j\geq  \sum_{i\neq j} \lambda _i\lambda _j; \;\;  \max \{d_i\}\geq \lambda _2.
$$

The S-NIEP and the R-NIEP are equivalent for $n\leq 4$, and remain unsolved for $n\geq 5$. 
Fiedler \cite[Theorem 4.8]{F} solved the S-NIEP in the case $n=3$ with fixed diagonal entries: 
$\{d_1,d_2, d_3\}$, with $d_1\geq d_2\geq d_3\geq 0$, is the diagonal of a 3-by-3 symmetric nonnegative matrix with spectra $\{\lambda _1, \lambda _2, \lambda _3\}$, where $\lambda _1\geq \lambda _2\geq \lambda _3$,
if and only if  $\{\lambda _1, \lambda _2, \lambda _3\}$ majorizes $\{d_1, d_2, d_3\}$ and $d_1\geq \lambda _2$.

Johnson-Laffey-Loewy \cite{JLL} showed that the R-NIEP and the S-NIEP are different, 
and Egleston-Lenker-Narayan \cite{ELN} proved that they are different for $n\geq 5$.

The S-NIEP for $n=5$ has been widely studied \cite{LoMc, McNe, ELN}, but not fully resolved. 
It is common to study it considering the number of positive eigenvalues. 
When there are 1, 4 or 5 positive eigenvalues the answer for the S-NIEP is straightforward. 
Recently, Johnson-Mariju\'an-Pisonero \cite{JMP2} resolved all cases with 2 positive eigenvalues, 
and they give a method, based upon the eigenvalue interlacing inequalities for symmetric matrices, 
to rule out many unresolved spectra with 3 positive eigenvalues. 
In particular, this method shows that the nonnegative realizable spectrum $\{6, 3, 3, -5, -5\}$ is not symmetrically realizable.

Also recently, Loewy-Spector \cite[Theorem 4]{LoSp} characterize the case $n=5$ in a particular case:
$\Lambda =\{\lambda_{1}, \lambda_{2}, \dots , \lambda_{5}\}$, where $\lambda_{1}\geq \lambda_{2}\geq \dots \geq \lambda_{5}$ and $2s_1(\Lambda )\geq  \lambda _1$, is (S-NIEP) realizable if and only if (Perron), 
$\lambda _2+\lambda _5\leq \mbox{Tr }(\Lambda ), \; \lambda _3\leq \mbox{Tr }(\Lambda ).$ This last condition is implied by the constraint $2s_1(\Lambda )\geq  \lambda _1$.

The trace 0 NIEP has also been extensively studied. Reams \cite{Re1} solved the case $n=4$: $\{\lambda_{1},\lambda_{2},\lambda_{3},\lambda_{4}\}$ is trace 0 (NIEP)-realizable if and only if $s_{1}=0, \;s_{2}, s_{3}\geq 0$ and $s_{2}^{2}\leq 4s_{4}$. The case $n=5$ was first studied by Reams \cite{Re1} and he gave a sufficient condition. The case $n=5$ was finally solved by Laffey-Meehan \cite{LM2}:  
$\{\lambda_{1},\dots ,\lambda_{5}\}$ is trace 0 (NIEP)-realizable if and only if
$s_{1}=0, \;s_{2}, s_{3}\geq 0,\; s_{2}^{2}\leq 4s_{4}$ and $12s_5+5s_3\sqrt{4s_4-s_2^2}\geq 5s_2s_3$.

Torre-Mayo {\it et al.} \cite{TAAMP} generalize these solutions in terms of the coefficients of the characteristic polynomial:
the polynomial $x^n+k_px^{n-p}+\cdots +k_{n-1}x+k_n$, with $2\leq p\leq n\leq 2p+1$, is (NIEP)-realizable if and only if 

\vskip0.3cm
\noindent 
$ k_p,\ldots ,k_{2p-1}\leq 0; \;k_{2p}\leq\frac{k_p^2}{4}; \; k_{2p+1}\leq\left\{
\begin{array}{ll}
	k_pk_{p+1}& \mbox{if }\;k_{2p}\leq 0,\\
	\vspace*{-.3cm} & \\
	k_{p+1}\Big(\frac{k_p}{2}-\sqrt{\frac{k_p^2}{4}-k_{2p}}\Big)&\mbox{if }\; k_{2p}>0.\\
\end{array}\right.$

%\noindent Moreover, when {\it i)} and {\it ii)} hold, $P(x)$ is EBL realizable.

\vskip0.3cm
 Spector \cite{Sp} characterized trace 0 S-NIEP realizability for $n=5$ by the conditions  $\lambda _2+\lambda _5\leq 0$ and $s_3\geq 0$.

%For all the NIEPS $n=1,2,3,4$. $n=5$ lots of specifics, importance of $\lambda_3$. Many recent results.

%----------------------------------------------------------------------------------------------------------------------------------------------------------------------------------------------------------------------------
\section{Sufficient Conditions} \label{Se5}
The first known sufficient condition for the NIEP, that in fact is for the R-NIEP, was announced by Sule\v{\i}manova \cite{Su} in 1949 and proved by Perfect \cite{Pe2} in 1953:
\begin{eqnarray*}
\left.\begin{array}{r}
\Lambda=\{\lambda_1,\dots,\lambda_n\} \mbox{ real, }\lambda_1\geq |\lambda| \mbox{ for } \lambda\in \Lambda  \\
\vspace*{-.3cm}\\
\mbox{ and } \lambda_1+\sum_{\lambda_i<0}\lambda_i\geq 0 
\end{array}\right\}
\Longrightarrow\; \Lambda\mbox{ is realizable}.
\end{eqnarray*}
Several other proofs have been given, {\it e.g.} \cite{Pa}. There are several sufficient conditions for the R-NIEP that are checkable in a straightforward way, that is, one only needs to check a few algebraic inequalities, perhaps after ordering the spectrum. The authors of them are: Ciarlet in 1968 \cite{Ci}, Kellogg in 1971 \cite{Ke}, Salzmann in 1972 \cite{Sa} and Fiedler in 1974 \cite{F}. Borobia in 1995 \cite{Bo1} extended  Kellogg's condition by  grouping negative eigenvalues.

Other sufficient conditions for the R-NIEP involve partitions of the spectra considered, such as an immediate piece-wise extension of  the Sule\v{\i}manova condition. We will name this condition {\bf  Sule\v{\i}manova-Perfect} \cite{Su,Pe2}:
\begin{eqnarray*}
\left.\begin{array}{r}
\Lambda=\{\lambda_1,\lambda_{11},\dots,\lambda_{1t_1},\dots,\lambda_r,\lambda_{r1},\dots,\lambda_{rt_r}\} \mbox{ real,}\\
\vspace*{-.3cm}\\
\lambda_1\geq |\lambda| \mbox{ for } \lambda\in \Lambda \\
\vspace*{-.3cm}\\
 \mbox{ and } \lambda_j+\sum\limits_{\lambda_{ji}<0}\lambda_{ji}\geq 0  \mbox{ for }j=1,\dots,r
\end{array}
 \right\}\Longrightarrow\; \Lambda\mbox{ is realizable}.
\end{eqnarray*}
Other sufficient conditions of this type are more elaborate and we will name them by their authors: {\bf Perfect 1} in 1953 \cite{Pe2}, Soto 2 in 2003 \cite{So1}  or its extension Soto $p$ in 2013 \cite{So5}. 

Some of the sufficient conditions besides partitions involve the knowledge of the diagonal entries of a realization of part of the spectrum. The first condition of this type is due to Perfect in 1955 \cite{Pe3}:
\begin{eqnarray*}
\left.\begin{array}{r}
\Lambda =\{\lambda _1, \dots, 
\lambda _r\}\cup\{\lambda _{11}, \dots,\lambda _{1t_1}\}\cup \dots\cup
\{\lambda _{r1},\dots,\lambda _{rt_r}\}\\
\vspace*{-.3cm}\\
\{\lambda _1, \dots, 
\lambda _r\} \mbox{ the spectrum of a nonnegative}\\
\vspace*{-.3cm}\\
\mbox{matrix with diagonal }d_1,\dots, d_{r},\\
\vspace*{-.3cm}\\
\lambda _{ji}\leq 0\mbox{  for } j=1,\dots,r \mbox{ and } i=1,\dots,t_j,\\
\vspace*{-.3cm}\\
\lambda _1 \geq |\lambda| \mbox{ for } \lambda \in \Lambda\;,\sum\limits_{\lambda \in \Lambda }\lambda \geq 0\;,\\
\vspace*{-.3cm}\\
\mbox{and }
d_j+\sum\limits_{1\leq i\leq t_j}\lambda _{ji}\geq 0\;,\; j=1,\dots, r
\end{array}
\right\} \Rightarrow \Lambda \mbox{ is realizable}.
\end{eqnarray*}
When $\lambda_j\geq 0$, for $j=1,\dots,r$, we call this condition {\bf Perfect 2$^+$} (see \cite{MPS}). There are two equivalent conditions that extend this condition: Soto-Rojo in 2006 \cite{SoRo} and Soto-Rojo-Manzaneda in 2011 \cite{SoRoMa}.

Other well-known sufficient conditions manipulate certain spectra to get a new realizable spectrum: Guo in 1997 \cite[Theorems 2.1 and 3.1]{G2} or $C$-realizability \cite{BMS2} that we name the {\bf game} condition \cite{MaPi2014}.

In order to construct a map of sufficient conditions for the R-NIEP, Mariju\'an-Pisonero-Soto compared these conditions and established inclusion relations or independence relations between them, \cite{MPS,MaPi2014}:

\vspace*{.15cm}

\hspace*{-.5cm} \begin{minipage}{.45\hsize}
{\pdfimage width \hsize {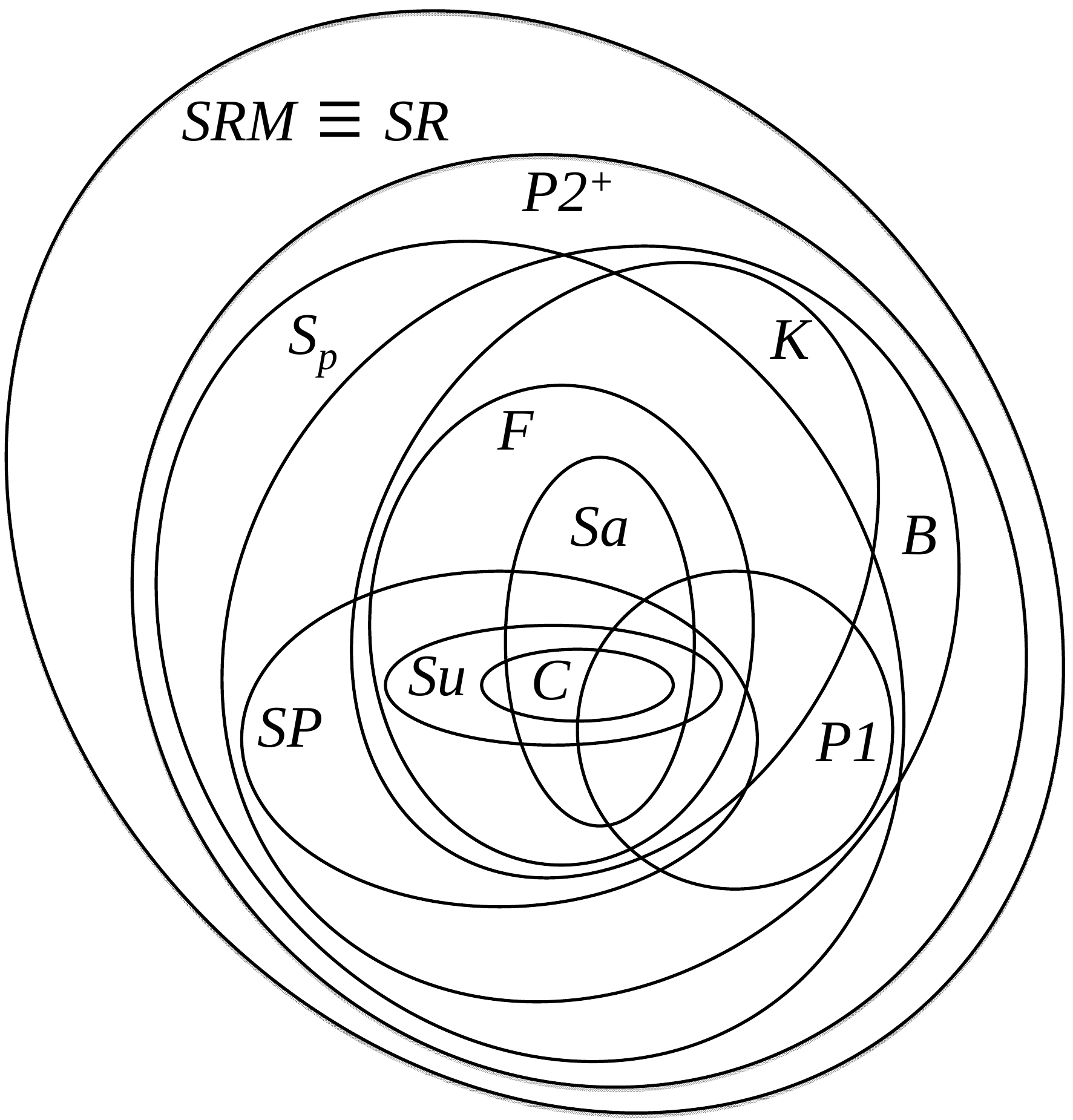}}
\end{minipage}\hspace*{.5cm}\begin{minipage}{0.5\hsize} 
\begin{itemize}
\item [\cite{Bo1}]  Borobia $= B$
%\vspace*{-.2cm}
\item [\cite {Ci}] Ciarlet $= C$
%\vspace*{-.2cm}
\item [\cite{F}] Fiedler  $= F$
%\vspace*{-.2cm}
\item [\cite{Ke}] Kellogg $= K$
%\vspace*{-.2cm}
\item  [\cite{Pe2}] Perfect 1 $= P1$
%\vspace*{-.2cm}
\item [\cite{Pe3}] Perfect 2$^+\;= P2^+$
%\vspace*{-.2cm}
\item [\cite{Sa}] Salzmann $= Sa$
%\vspace*{-.2cm}
\item [\cite{So5}] Soto p $= S_p$
%\vspace*{-.2cm}
\item [\cite{SoRo}] Soto-Rojo $= SR$
%\vspace*{-.2cm}
\item [\cite{SoRoMa}] Soto-Rojo-Manzaneda $= SRM$
%\vspace*{-.2cm}
\item [\cite{Su}] Sule\v{\i}manova $= Su$
%\vspace*{-.2cm}
\item [\cite{Pe2}] Sule\v{\i}manova-Perfect $= SP$
\end{itemize}
\end{minipage} 

\vspace*{.2cm}

\hspace*{-.7cm}\begin{minipage}{.45\hsize}
\setlength{\unitlength}{1cm}
\begin{picture}(4,4)
{\pdfimage width 1.1\hsize {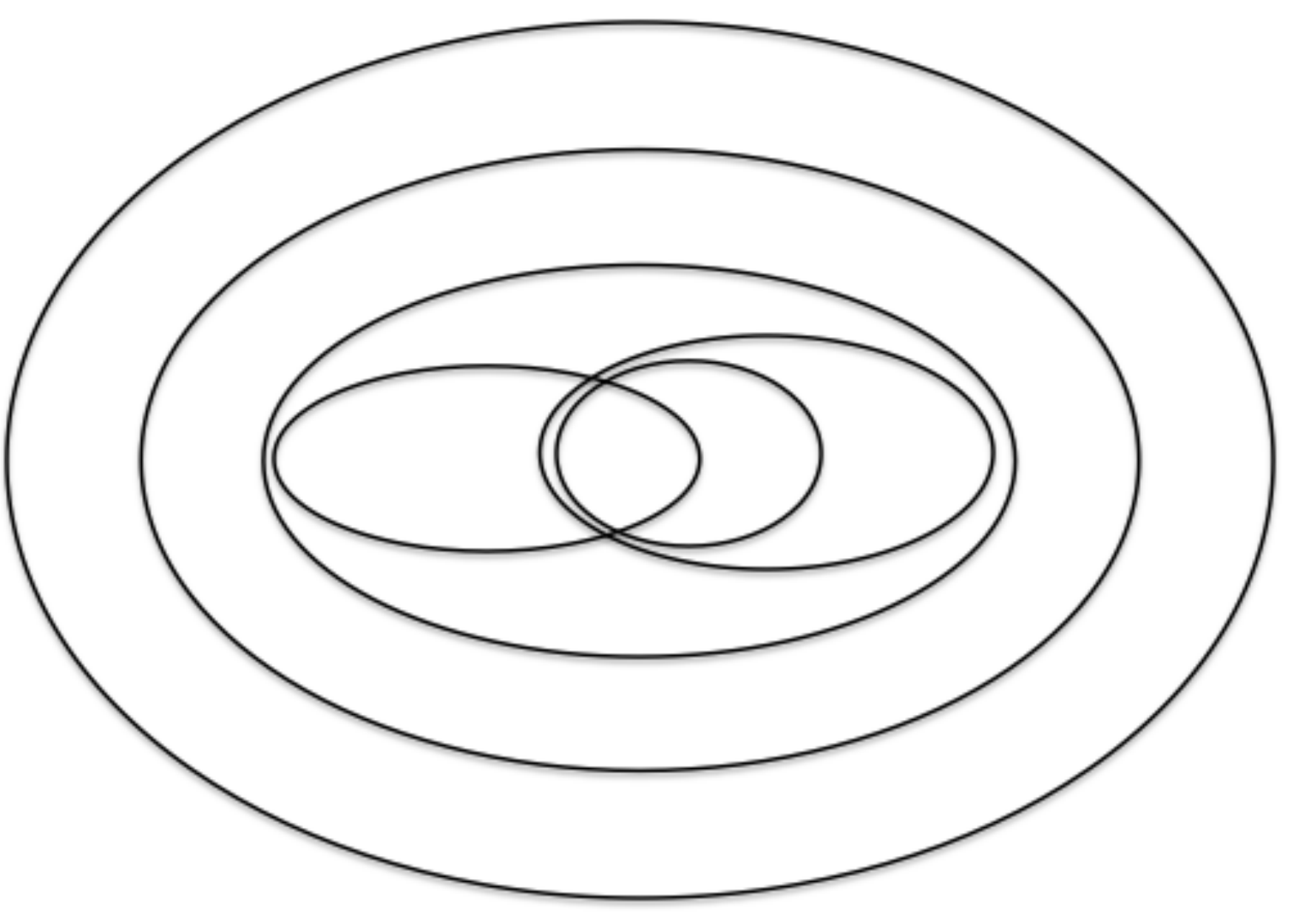}}
\put(-3.9,2){{$B$}}
\put(-2.7,2){{ $S_p$}}
\put(-2.1,2){{$S_{p+1}$}}
\put(-3.7,2.6){{ Sotos}}
\put(-3.3,3.1){{ game}}
\put(-3.1,3.65){{ $P2^+$}}
\put(-1.4,2){  $\overset {\mbox{\bf?}}{\longleftrightarrow}$}
\end{picture}
\end{minipage}\hspace*{.9cm}\begin{minipage}{0.5\hsize} 
\begin{itemize}
\item [\cite{Bo1}]   Borobia $= B$
%\vspace*{-.2cm}
\item [\cite{Pe3}] Perfect 2$^+\;= P2^+$
%\vspace*{-.2cm}
\item [\cite{So5}] Soto p $= S_p$
%\vspace*{-.2cm}
\item [\cite{So5}] Soto p + 1 $= S_{p+1}$
%\vspace*{-.2cm}
\item [\cite{So5}] Sotos = $\bigcup\limits_{p\geq 2} Soto\; p$
\end{itemize}
\end{minipage} 

\vspace*{.15cm}

The first known sufficient condition for the S-NIEP is due to Perfect-Mirsky in 1965 \cite{PeMi} for doubly stochastic matrices, and Fiedler in 1974 \cite{F} gave the first one for symmetric nonnegative matrices. Several sufficient conditions which were first obtained for the R-NIEP have later been shown to be valid also for the S-NIEP as well. Fiedler \cite{F}, Radwan \cite{Ra} and Soto \cite{So3} showed, respectively, that Kellogg \cite{Ke}, Borobia \cite{Bo1} and Soto 2 \cite{So1} are also symmetric sufficient conditions. 

Soules in 1983 \cite{Soules} gave two constructive sufficient conditions for symmetric realization. The inequalities that appear in these conditions are obtained by requiring the diagonal entries of the matrix $R\mbox{diag}(\lambda_1,\dots,\lambda_n)R^T$ to be nonnegative, in which $R$ is an orthogonal matrix with a certain pattern. For a particular $R$, this condition is the Perfect-Mirsky condition.

Soto-Rojo-Moro-Borobia gave in 2007 \cite{SRMB} a symmetric version of the Soto-Rojo condition. Soto $p$ \cite{So5} and Soto-Rojo-Manzaneda \cite{SoRoMa} have also symmetric versions. Laffey-\v{S}migoc in 2007 \cite{LS2} gave the symmetric realizability of a spectrum by manipulating two spectra.

Again in order to construct a map of sufficient conditions for the S-NIEP, Mariju\'an-Pisonero-Soto \cite{MPS1} compared these conditions and established inclusion relations or independence relations between them:

\hspace*{-.8cm}\begin{minipage}{.55\hsize}
{\pdfimage width \hsize {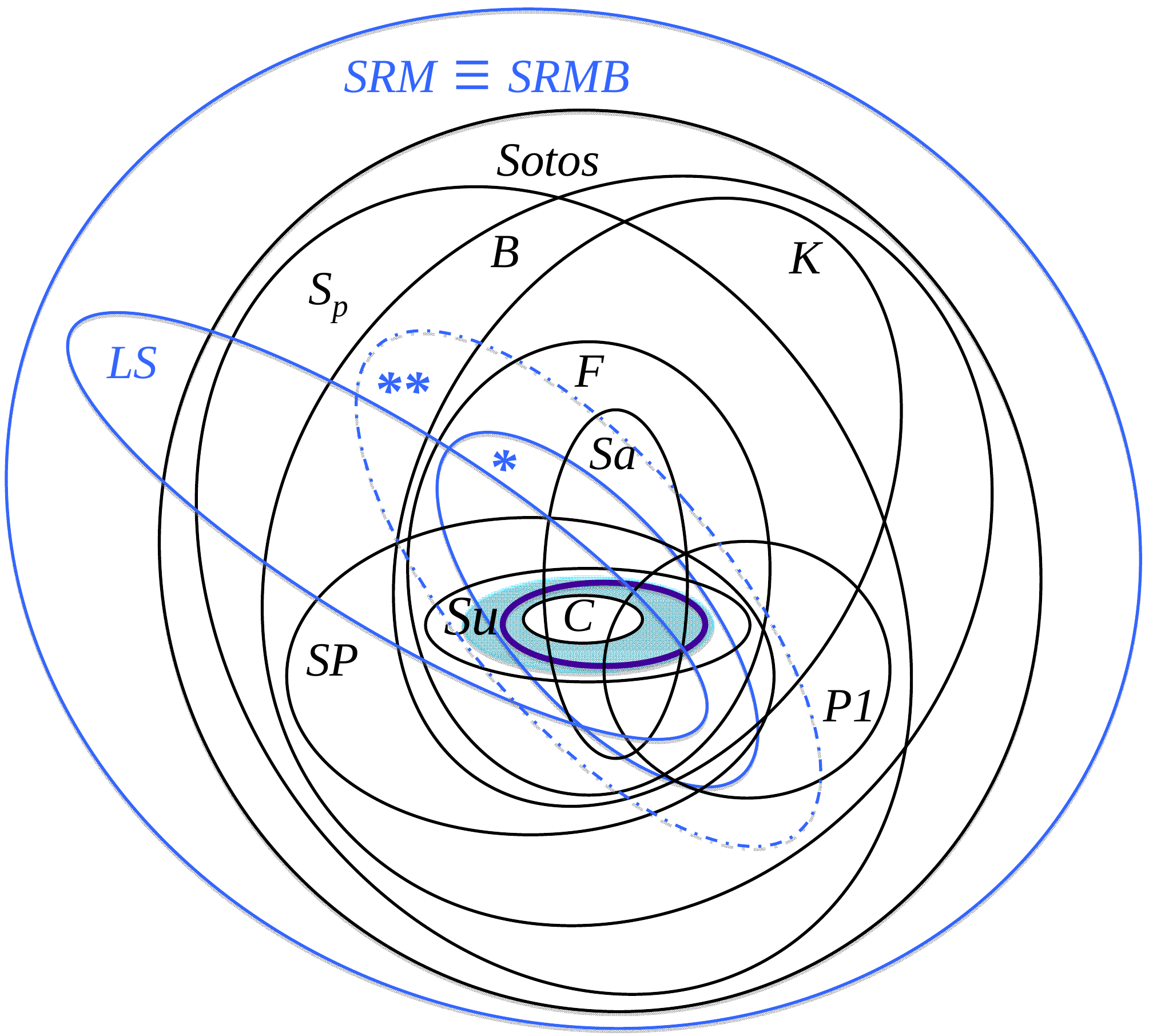}}
\end{minipage}\hspace*{.3cm}
\begin{minipage}{0.45\hsize} 
\begin{itemize}
\item [\cite{Bo1}]  Borobia $= B$
%\vspace*{-.2cm}
\item [\cite {Ci}] Ciarlet $= C$
%\vspace*{-.2cm}
\item [\cite{F}] {\color{blue}Fiedler $= $} $ F$
%\vspace*{-.2cm}
\item [\cite{Ke}] Kellogg $= K$
%\vspace*{-.2cm}
\item [\cite{LS2}] {\color{blue} Laffey-\v{S}migoc $= LS$}
%\vspace*{-.2cm}
\item  [\cite{Pe2}] Perfect 1 $= P1$
%\vspace*{-.2cm}
\item  [\cite{PeMi}] {\color{blue} Perfect-Mirsky $=$}  {\tiny ${\color{morado} \blacksquare \hspace*{-.1cm}\blacksquare\hspace*{-.1cm}\blacksquare \hspace*{-.1cm}\blacksquare\hspace*{-.1cm}\blacksquare\hspace*{-.1cm}\blacksquare\hspace*{-.1cm}\blacksquare}$}
%\vspace*{-.2cm}
\item [\cite{Sa}] Salzmann $= Sa$
%\vspace*{-.2cm}
\item [\cite{So5}] Soto p $= S_p$
%\vspace*{-.2cm}
\item [\cite{So5}] Sotos = $\bigcup\limits_{p\geq 2} Soto\; p$
%\vspace*{-.2cm}
\item [\cite{SoRoMa}]  {\color{blue}Soto-Rojo-Manzaneda $= SRM$}
%\vspace*{-.2cm}
\item [\cite{SRMB}]  {\color{blue}Soto-Rojo-Moro-Borobia $= SRMB$}
%\vspace*{-.2cm}
\item [\cite{Soules}] {\color{blue} Soules 1 $=$} {\Large ${\color{azulclaro} \blacksquare}$}
%\vspace*{-.2cm}
\item [\cite{Soules}] {\color{blue} Soules 2 $=$  {\Large $**$}}
%\vspace*{-.2cm}
\item [\cite{Soules}] {\color{blue} Soules 2 corollary $=$  {\Large $*$}}
%\vspace*{-.2cm}
\item [\cite{Su}] Sule\v{\i}manova $= Su$
%\vspace*{-.2cm}
\item [\cite{Pe2}]  Sule\v{\i}manova-Perfect $= SP$
\end{itemize}
\end{minipage} 

\vspace*{.5cm}

\noindent The discontinuous line for Soules 2 in the map means that we only conjecture this position for this sufficient condition.

Recently, Ellard-\v{S}migoc \cite{ElSm2} have modified the Laffey-\v{S}migoc  condition and the Soules 2 condition and they have proved the following equivalence:
\begin{eqnarray*}
\mbox{Soules 2 modified } \Longleftrightarrow \mbox{ Laffey-\v{S}migoc modified } \Longleftrightarrow \mbox{ game } \Longleftrightarrow \mbox{ Sotos}.
\end{eqnarray*}
This implies the symmetric realizably of the game condition and resolves the question on the RNIEP diagram about whether there exists something between Sotos and game.

There are also some sufficient conditions for the NIEP. Guo in 1997 \cite[Theorems 2.1 and 3.1]{G2} gave results about changing a realizable spectrum to obtain a realizable spectrum. \v{S}migoc in 2004 \cite[Theorems 10 and 14]{Sm1} gave other such results. Finally, let us mention some generalizations of the Sule\v{\i}manova condition: Laffey-\v{S}migoc in 2006 \cite[Theorems 1 and 3]{LS1} and Borobia-Moro-Soto in 2008 \cite[Theorem 3.3]{BMS2}.
%%%%%%%%%%%%%%%%%%%%%%%%%%%%%%%%%%%%%%%%%%%

%\underline{NIEP}: Generalization of Sule\v{\i}manova \cite[Theorem 3.3]{BMS2}, Guo 97 \cite[Theorems 2.1 and 3.1]{G2}, \v{S}migoc 2004 \cite[Theorems 10 and 14]{Sm1}, Laffey-\v{S}migoc 2006 \cite[Theorems 1 and 3]{LS1} (to be considered in Section 6)

%\underline{RNIEP}:

%Straight conditions: write Sule\v{\i}manova (permutative proof) as the first condition and mention Ciarlet, Salzmann and Fiedler in a map that relates them.

%Sule\v{\i}manova (permutative proof), majorization, variants. Soules, Fiedler, Kellogg, etc. Interrelations.

%----------------------------------------------------------------------------------------------------------------------------------------------------------------------------------------------------------------------------
\section{Embedding Spectra, by Adding 0's, to Achieve Realizability} \label{Se6}
Consider now spectra that meet the most basic necessary conditions  for NIEP-realizability:
(Perron), (Trace) and  (Reality). Even when these conditions are strictly met, the proposed
spectrum, $\sigma $, need not be realizable (for example $\{3,3,-\sqrt{3}\pm i\}$). This raises
 the natural question of whether the spectrum may be embedded
in a  larger one, that is realizable, by appending some additional
eigenvalues, {\it e.g.} $\sigma \hookrightarrow \sigma \cup \tau$.
If we are too liberal about what eigenvalues may be appended, this
question becomes trivial. For example, realizability may always be
achieved by appending a single, sufficiently large, positive
eigenvalue; the Perron and trace conditions may be arbitrarily
improved. Thus, some condition must be placed upon the appended
eigenvalues. A natural one is that only 0 eigenvalues may be
appended; now, the Perron and trace conditions are not enhanced -
but, the dimension is increased. (Intermediate restrictions seem
not yet to have been considered.) The increase in dimension does
improve the possibility of meeting the JLL conditions.

 The notion of appending 0's to ''repair'' a nonrealizable spectrum
may and has often been viewed another way: what collections of
complex numbers occur as the ''nonzero part'' of the spectrum of an entry-wise nonnegative matrix. It turns out that this
question  is quite different from and more tractable than the
classical  NIEP.

The first to show that appending 0's
can help was \cite{J2}, in which it was shown that the
spectrum
\begin{eqnarray*}
1,\sqrt{\frac{3}{8}}i,-\sqrt{\frac{3}{8}}i
\end{eqnarray*}
is not realizable in dimension 3 (because the single eigenvalue
conditions \cite{Ka} are not met, or the JLL conditions are
not  met), but the spectrum
\begin{eqnarray*}
1,\sqrt{\frac{3}{8}}i,-\sqrt{\frac{3}{8}}i,0
\end{eqnarray*}
is realizable in dimension 4 by the matrix
\begin{eqnarray*}
\left( \begin{array}{cccc}
 1/4&0&3/4&0\\
 1/4&1/4&0&1/2\\
 0&3/4&1/4&0\\
 7/24&0&11/24&1/4
 \end{array}\right).
\end{eqnarray*}

 Later in \cite{BH}, it was shown, remarkably, that if a
spectrum  (with no 0's) meets the basic necessary conditions
strictly, then  it is the nonzero part of the spectrum of a
nonnegative matrix.  Sufficiently many 0's may be appended to
achieve realizability. Of  course, the number of 0's that need
be added may be very large  (not uniformly bounded in terms of
$n$) because of JLL.  Interestingly, the newish method of
symbolic dynamics was used in  an import way, though many
easily proven matricial lemmas were  needed as well.
Specifically the result is

\begin{teo}
 The list of nonzero complex numbers $\Lambda
=\{\lambda_1,\lambda_2,\dots,\lambda_n\}$ is the nonzero spectrum
 of a primitive matrix if and only if
\begin{enumerate}
\item {\bf strict Perron condition}: $\Lambda$ contains a positive
 eigenvalue of multiplicity one that is greater in absolute
value  than all other $\lambda_i's$;
\item {\bf reality
condition}:  $p(t)=\prod\limits_{i=1}^{n}(t-\lambda_i)$ has real
coefficients;

and

\item {\bf extended trace condition}: $\sum\limits_{i=1}^{n}\lambda_i^{mk}\geq 0$ for all $k$ and
$\sum\limits_{i=1}^{n}\lambda_i^{k}>0$ implies
$\sum\limits_{i=1}^{n}\lambda_i^{mk}>0$ for all $m$.
\end{enumerate}
\end{teo}

The necessity of
these conditions is easily verified (the last one via JLL, for
example), and the interesting point is that necessary conditions
 become sufficient when primitivity is the goal and the
dimension  may be arbitrarily increased. When the Perron
condition is not  strict it may not be possible to ''save'' a
spectrum meeting  obvious necessary conditions. The familiar
example $3,3,-2,-2,-2$  is not only not realizable, but is never
the nonzero part of the  spectrum of a nonnegative matrix.
However, both  $3+\epsilon,3,-2,-2,-2$ and
$3+\epsilon,3-\epsilon,-2,-2,-2$ are  both the nonzero parts of
the spectra of primitive matrices for  arbitrarily small
$\epsilon >0$ (though they are not realizable  for $\epsilon$
small).

% (should we have a small module about $3,3,-2,-2,-2$ and  $3+\epsilon,3,-2,-2,-2$ for the R-NIEP and S-NIEP here or  elsewhere?)

 More recently
there have been further developments about  realizability after
appending 0's to a spectrum. In \cite{La3}  there is matricial
proof of key results from \cite{BH}, which is  much more
explicit. Though the number of 0's needed to make a  spectrum
realizable may be very large (and not easy to estimate  from
\cite{BH}), estimates have recently been given under some circumstances by bringing $s_2(\Lambda)$,
 as well as the trace, into play \cite{LS1}.

%History, CJ exp, precise results/questions, recent improvements (unbounded number, relationship with J-L-L).

\section{The Graph NIEP} \label{Se7}
Not surprisingly, realizable spectra that are, in some way
extremal, are often realizable by nonnegative matrices with many 0
entries, as are many non-extremal spectra. This raises the
question that, if we fix the 0-pattern of a nonnegative matrix,
how is the NIEP restricted, {\it i.e.} which realizable spectra
occur? A natural way to describe a particular 0-pattern is via a
graph, which could be directed or undirected. For a particular
graph $G$ on $n$ vertices, consider the set of nonnegative
$n$-by-$n$ matrices $N(G)$ for which $A=(a_{ij})$ satisfies
$a_{ij}>0$ if and only if $(i,j)$ is a directed edge, $i\neq j$, of $G$.
(Here, we consider an undirected edge $\{i,j\}$, if $G$ is
undirected, to consist of two directed edges $(i,j)$ and $(j,i)$.)
For simplicity we consider graphs without loops, and no
restriction is placed by the graph upon the diagonal entries,
other than nonnegativity. The G-NIEP then just asks which spectra
occur among matrices in $N(G)$? (if we wish to emphasize the
dimension, which we generally take to be implicit, we may write
$N_n(G)$.) Of course the NIEP is just the union of the solutions to
the G-NIEP's, over all directed graphs $G$. The same is true for
various variations upon the NIEP. For example, the solution to the
S-NIEP is just the solution to the G-NIEP, restricted to symmetric
matrices, over all undirected graphs. Thus far, the G-NIEP has
been considered only for undirected graphs $G$, symmetric matrices
and real eigenvalues. This seems a fertile area for future work.

A prototype of the G-NIEP, though not presented in graph terms,
appeared in \cite{FriMe}, in which tridiagonal matrices were
considered. Notice that the spectrum of a nonnegative tridiagonal
matrix is not only necessarily real but also the spectrum of a
symmetric tridiagonal nonnegative matrix. Also, tridiagonal
matrices are the case in which $G$ is a path (and an edge only
requires nonnegativity of the entry); off-diagonal 0 entries are
important. We say that a matrix $A$ is {\bf subordinate} to a
graph $G$ if $G(A)$ has the same vertices as $G$ and the edges of
$G(A)$ are contained among those of $G$; equivalently, for
$A=(a_{ij})$, $a_{ij}\neq 0$ implies $\{i,j\}$ is an edge of $G$. The main result of \cite{FriMe} is for the path $P$ on $n$ vertices. Then, for $\lambda_1\geq \lambda_2\geq \cdots \lambda_n$, there is an $n$-by-$n$ nonnegative matrix $A^T=A$ subordinate to $P$ and with eigenvalues $\lambda_1,\lambda_2,\dots , \lambda_n$ if and only if
\begin{eqnarray*}
\lambda_i+\lambda_{n-i+1} \geq 0
\end{eqnarray*}
$i=1,2,\dots,n$. Additional conditions, when the graph is precisely $P$, are also discussed, but the results are incomplete.

In \cite{ LeJo,J2}, the observation of \cite{FriMe} is dramatically generalized. Recall that a path is a tree and that all trees are bipartite (but many non-trees are bipartite as well). All bipartite graphs are considered in \cite{ LeJo}. The main result is that $\lambda_1\geq \cdots \geq  \lambda_n$ are the eigenvalues of an $n$-by-$n$ nonnegative symmetric matrix $A$ subordinate to a given bipartite graph $G$ on $n$ vertices if and only if
\[ \lambda_1+\lambda_n  \geq  0\]
\[\lambda_2+\lambda_{n-1} \geq 0\]
\[\vdots\]
\[\lambda_m+\lambda_{n-m+1} \geq 0\]
\[\lambda_{m+1}\geq 0\]
\[\vdots\]
\[\lambda_{n-m}\geq 0\]
in which $m$ is the matching number of $G$. (Note that $\lfloor \frac{n}{2}\rfloor$ is the matching number of a path on $n$ vertices, so that the result of \cite{FriMe} is a special case.) Further observations about the S-NIEP for matrices subordinate to a given graph $G$ are also made.

%Friedland-Melkman for tridiagonal, L-D, J for bipartite, etc

\section{Perron Similarities} \label{Se8}
%Hadamard matrices, use, sample cases, $n=4$, Relative gain array. RHC matrices and cone containment, etc

If $S$ is an invertible matrix and there is a real, diagonal, nonscalar matrix $D$ with $S D S^{-1} \geq 0$, then $S$ is called a \emph{Perron similarity}. Perron similarities were introduced by Johnson and Paparella to study the diagonalizable R-NIEP and the S-NIEP. 

Perron similarities were characterized in several ways in \cite{JoPa1} and it was shown that $\mathcal{C}(S):= \{ x \in \R^n: S D_x S^{-1} \geq 0 \}$ is a \emph{polyhedral cone}, {\it i.e.}, a convex cone in $\R^n$ with finitely-many extremals. This was called the \emph{(Perron) spectracone of $S$}, and a certain cross-section, a polytope called the \emph{(Perron) spectratope of $S$}, was also discussed. 

These polyhedral sets were used to verify the known necessary and sufficient conditions for the R-NIEP and the S-NIEP for orders up to four. For orders $1$, $2$, and $4$, it is shown that a finite number of Perron similarities are required to cover the realizable region, whereas when $n=3$, it is shown, via the \emph{relative gain array} (see, {\it e.g.}, \cite{HoJo2}), that an uncountable number of Perron similarities is required to cover the realizable region. 

For every $n \geq 1$, the spectracone and spectratope of $H_n$ were characterized, where $H_n$ denotes the \emph{canonical Hadamard matrix of order $2n$}. More specifically, the spectracone of $H_n$ is the conical hull of its rows and the spectratope of $H_n$ is the convex hull of its rows. The spectratope of a general, normalized Hadamard matrix was used to give a constructive proof, for \emph{Hadamard orders}, of a result by Fiedler \cite{F}, that every \emph{Sule\v{\i}manova spectrum} is the spectrum of a symmetric nonnegative matrix (in fact, this result was strengthened to show that the constructed realizing matrix is also \emph{doubly stochastic}). It is still an open problem to find a constructive proof that every Sule\v{\i}manova spectrum is SNIEP-realizable. A constructive proof of the Boyle-Handelman theorem for Sule\v{\i}manova~spectrum: augmenting any such spectrum with zeros up to a Hadamard order yields a spectrum realizable by a nonnegative matrix that is symmetric and doubly stochastic.

\section{Jordan structure and the NIEP's} \label{Se9}
When there are repeated eigenvalues in a proposed spectrum, a natural question is whether Jordan structure for the repeated eigenvalues can play a role in realizability. The J-NIEP asks which particular Jordan canonical forms occur for $n$-by-$n$ nonnegative matrices and the D-NIEP is the special case of which spectra occur among diagonalizable $n$-by-$n$ nonnegative matrices. Of course, the D-NIEP and the NIEP are the same for spectra with distinct eigenvalues. It is also a simple exercise that any spectrum that is D-NIEP realizable by a positive matrix is also J-NIEP realizable for any Jordan canonical form possible for its eigenvalues. However, whether D-NIEP realizability always implies J-NIEP realizability for any possible Jordan form is unclear. This is so for Sule\v{\i}manova spectra \cite{CcSo1, SDNS, DiSo}.

There are spectra that are realizable (in fact, R-NIEP realizable), but are not diagonalizably realizable. The smallest dimension in which this occurs is 5. For $n\leq 4$, any realizable spectrum with (non-real) complex eigenvalues is diagonalizably realizable (the only difficulty could come from a multiple Perron root and that is easily handled by a reducibility argument). For $n\leq 4$, R-NIEP and S-NIEP realizability are the same, which settles the matter. However, for $n=5$, the spectrum
\[3+t,3-t,-2,-2,-2\]
is realizable for $t>(16\sqrt{6})^{1/2}-39\approx 0.437...$, \cite{LM2}. However, it is diagonalizably realizable iff $t\geq 1$ 
%\cite{..} Cronin-Laffey to appear 
(in which case it is also symmetrically realizable). Thus, for $0.437...<t<1$, this 5-spectrum is realizable, but not diagonalizable so. This suggests that this phenomenon is fairly common.

%%%%%%%%%%%%%%%%%%%%%%%%%%%%%

\end{document}